# Gambler's Ruin Random Walks and Brownian Motions in Reserves Modeling. Application to Pensions Funds Sustainability

**Manuel Alberto M. Ferreira** [1] **and José António Filipe** [2,]

[1] Department of Mathematics, ISTA—School of Technology and Architecture, Iscte - Instituto Universitário de Lisboa; Information Sciences, Technologies and Architecture Research Center (ISTAR-IUL); Business Research Unit-IUL (BRU-IUL), 1649-026 Lisbon, Portugal; manuel.ferreira@iscte-iul.pt

[2] Department of Mathematics, ISTA—School of Technology and Architecture, Iscte - Instituto Universitário de Lisboa; Information Sciences, Technologies and Architecture Research Center (ISTAR-IUL); Business Research Unit-IUL (BRU-IUL), 1649-026 Lisbon, Portugal; jose.filipe@iscte-iul.pt

**Abstract:**

We used the random walk to model the problem of reserves. The classic case of a stochastic process is the example of random walks, which are used to study a set of phenomena and, particularly, as in this article, models of reserves evolution. Random walks also allow the construction of significant complex systems and are also used as an instrument of analysis, being used in the sense of giving a theoretical characteristic to other types of systems. Our goal is primarily to study reserves to see how to ensure that pension funds are sustainable. This classic approach to the study of pension funds makes possible to draw interesting conclusions about the problem of reserves.

**1. Introduction**

The quality and the level of pensions is an important constituent in population's socio-economic status, see Iashina, Petukhov and Lazarev [1]. However, pension plan systems around the world have been facing many difficulties, as referred in Cristea and Thalassinos [2]. These authors show that the value of pension fund assets around the world was considerably reduced with the economic and financial crisis that began in 2008 and that caused a restructuring of pension fund assets and investments toward low-risk ones and often more directed to domestic investments, mainly in the OECD countries, see again Cristea and Thalassinos [2]. If, in the future, there is a lower pensions funds' capability for keeping a proper protection basis, that will make many population's segments with no sufficient income in old ages get at risk of poverty.

Pensions investments strategies are nowadays particularly complex and in the current days it is necessary to have an important understanding of this phenomenon, once pension assets are one of the largest pools of investment assets in the world, as highlighted in "Pensions: Global Issues, Perspectives and Challenges", approaching a set of specific problems and contexts, Webb [3] and Filipe [4].

As pension systems aim to protect economically and financially older people, allowing them a decent living, the pensions sustainability is crucial in this sense. It is necessary not forget that, as soon as aging society deepens as happens in many developed countries, it is expected that there may occur negative pension cash flows since contributions to the fund are lesser than the payments for pensions. This highlights the need for a very responsible investments strategy and a very qualified management policy.

To develop our study, the gambler's ruin statistical concept is applied in this paper, expressing the idea that gamblers play a negative expected value game that can eventually go broken.

Gambler's ruin suggests that growth events occur randomly, depending their survival on the stock of stored resources. In Gambler's ruin problem, reserves behave according to a simple random walk. This problem has been often presented in many works considering the stochastic processes theory in contexts under the frameworks of Markov Chains, Random Walks, Martingales or even others. Billingsley [5] or Feller [6] solve problems in this topic using the classic first step analysis to obtain a difference equation, which approach we use in our study. In another basis, Grimmett and Stirzaker [7] and Karlin and Taylor [8] get resolutions through the Martingales Theory and use it as an applications' example of the Martingales Stopping Time Theorem.

The study of the problem of the pension's funds sustainability, in the context it is developed, considers the gambler's ruin probability, i.e. the probability of the game reserves exhaustion. This theoretical basis brings an important insight over the problem that is a significant concern of governments' social politics, following some developments that have been studied in our previous works, as will be seen next. In this work we bring a new perspective and some complements to the previous approach, presenting a conception's repositioning and showing the usefulness of our methodology by bringing for instance important recommendations in this area and some orientations in terms of the established principles and future bases for funds management.

For national pensions funds plans, there is not an alternative way to a balanced pensions fund. Countries must give great attention to the long-term pension plans' financial sustainability, as many are doing, and it is this capability of getting a balanced long-term plan that will allow this systems' sustainability, Chen and Yang [9]. Besides, pensions' plans, whichever they are public or private, are ruled by the same needed principle of sustainability. In our paper, we provide tools for assessing the costs of policies followed to ensure the sustainability of these funds.

The paper is an updated approach of the work presented in SPMCS 2012, see Filipe, Ferreira, and Andrade [10]. We provide new proposals to develop the theoretical bases and to redefine some orientations concerning the consequent implications on pensions funds. After the presentation of "Materials and Methods", the "Gambler's Ruin Probability" is calculated, and then the "Fund Ruin Probability" is studied. Then follows a "Ruin Probability's Particularization" using normal distributions, the definition of Assets and Liability Management Politics the "Discussion" and ending the "References".

In Ferreira [11,12], Ferreira and Andrade [13], Ferreira, Andrade, and Filipe [14] and Andrade, Ferreira and Filipe [15] more applications of stochastic processes can be seen in problems of reserves and pensions funds.

**2. Materials and Methods**

Our method has its foundations in stochastic processes and in a set of concepts presented in the current section. These concepts allow us to reach results that show the importance of balanced pensions funds, supported on gambler's ruin modeling and on Brownian motion processes modeling.

We begin by considering a gambler with an initial capital of *x* monetary units, intending to play a sequence of games till the gambler's wealth reaches *k* monetary units. We suppose that *x* and *k* are integer numbers satisfying the conditions *x* > 0 and *k* > *x*. In each game, there are two possible situations: winning 1 monetary unit with probability *p* or losing 1 monetary unit with probability *q* = 1- *p*. A question is posed in terms of knowing which will be the probability the gambler gets in ruin before attaining his/her objective. This means that it is intended to know which the probability of losing the *x* monetary units is before adding gains of *k* – *x* monetary units.

Be $X_n, n = 1, 2, ...$ the outcome of the *n*th game. Obviously, the variables $X_1, X_2, ...$ are i.i.d. random variables, which common probability function is: $P(X_n = 1) = p, P(X_n = -1) = q = 1 - p$.

Consequently, the gambler's wealth, his/her reserves after the *n*th game represent the simple random walk as follows: $S_0 = x, S_n = S_{n-1} + X_n, n = 1, 2, ...$ .

We also consider a fund in which contributions/pensions received/paid, per time unit, are described as a sequence of random variables $\xi_1, \xi_2, ... \; (\eta_1, \eta_2, ...)$. Let us consider that $\xi_n (\eta_n)$ is the amount of the contributions/pensions that are received/paid by the fund during the $n^{th}$ time unit and consequently $X_n = \xi_n - \eta_n$ is the reserves variation in the fund at the $n^{th}$ time unit. If $X_1, X_2, ... X_n$ is a sequence of non-degenerated i.i.d. random variables, the stochastic process, representing the evolution of the fund reserves, since the value *x* till the amount $\tilde{S}_n$ after *n* time units, will be defined recursively as $\tilde{S}_0 = x, \tilde{S}_n = \tilde{S}_{n-1} + X_n$, with $n =, 1,2, ...$ . Such a process is a general random walk.

Finally, we assume that the assets value process of a pensions fund may be represented by the geometric Brownian motion process $A(t) = b \, e^{a+(\rho+\mu)t+\sigma B(t)}$ with $\mu < 0$ and $ab\rho + \mu\sigma > 0$, where $B(t)$ is a standard Brownian motion process. Suppose also that the fund liabilities value process performs such as the deterministic process $L(t) = be^{\rho t}$. Also consider the stochastic process $Y(t)$ obtained by the elementary transformation of $A(t)$, $Y(t) = \ln \frac{A(t)}{L(t)} = a + \mu t + \sigma B(t)$. It is a generalized Brownian motion process, starting at *a*, with drift μ and diffusion coefficient $\sigma^2$. Note also that the firs passage time of the assets process A(t) by the mobile barrier $T_n$, the liabilities process, is the first passage time of Y(t) by 0 with finite expected time. This framework is used to get an asset-liability management scheme of a pensions fund.

## 3. Results

### 3.1 Gamblers' Ruin Probability

Aiming to get the gambler's ruin probability, let us consider this probability as $\rho_k(x)$. It relates to the probability that $S_n = 0$ and $0 < S_i < k, i = 0, 1, ..., n - 1$ for $n = 1$ or $n = 2$ or .... If $\rho_k(x)$ is conditioned to the result of the first game, and considering the Total Probability Law, we get the following:

$$\rho_k(x) = p\rho_k(x + 1) + q\rho_k(x - 1), \qquad 0 < x < k \qquad (3.1.1).$$

If $0 \leq x \leq k$, the difference equation presented is easily solved considering the border conditions

$$p_k(0) = 1, \; p_k(k) = 0 \qquad (3.1.2).$$

So (3.1.1) can now be written as

$$p_k(x) - p_k(x-1) = \frac{p}{q}\left(p_k(x+1) - p_k(x)\right), \qquad 0 < x < k \qquad (3.1.3).$$

Considering $x = k - 1$ and also equation (3.1.2), $p_k(k-2) = p_k(k-1)\left(1 + \frac{p}{q}\right)$. Then, if $x = k - 2$, (3.1.3) become $p_k(k-3) = p_k(k-1)\left(1 + \frac{p}{q} + \left(\frac{p}{q}\right)^2\right)$. Continuing the process, the following general expression is achieved:

$$p_k(k-y) = p_k(k-1)\left(1 + \frac{p}{q} + \left(\frac{p}{q}\right)^2 + \cdots + \left(\frac{p}{q}\right)^{y-1}\right), 0 < y \leq k \qquad (3.1.4).$$

By means of (3.1.2) again, and considering (3.1.4), with $y = k$

$$p_k(k-1) = 1/\left(1 + \frac{p}{q} + \left(\frac{p}{q}\right)^2 + \cdots + \left(\frac{p}{q}\right)^{k-1}\right) \qquad (3.1.5).$$

Substituting (3.1.5) in (3.1.4) and accomplishing the variable change $y = k - x$ we reach finally the difference equation (3.1.1) solution with the border conditions (3.1.2):

$$p_k(k-1) = \begin{cases} \dfrac{1 - (p/q)^{k-x}}{1 - (p/q)^k}, & \text{if } p \neq \dfrac{1}{2} \\[2ex] \dfrac{k-x}{k}, & \text{if } p = \dfrac{1}{2} \end{cases} \qquad (3.1.6).$$

Be $N_a$ the random walk $S_n$ first passage time by $a$: $N_a = \min\{n \geq 0: S_n = a\}$. Now we can write $p_k(x) = P(N_0 < N_k | S_0 = x)$. It is appropriate to consider in (3.1.6) the limit as $k$ converges to $\infty$ to evaluate $p(x)$, the ruin probability of a gambler infinitely ambitious. In this context, considering the simple random walk $S_n$, $p(x) = P(N_0 < \infty | S_0 = x)$, after (3.1.6)

$$\rho(x) = \lim_{k \to \infty} \rho_k(x) = \begin{cases} (q/p)^x, & \text{if } p > \frac{1}{2} \\ 1, & \text{if } p \leq \frac{1}{2} \end{cases} \quad (3.1.7).$$

Being $\mu = E[X_n] = 2p - 1$, after (3.1.7), the ruin probability is 1 for the simple random walk at which the mean of the step is $\mu \leq 0 \Leftrightarrow p \leq \frac{1}{2}$.

## 3.2. Fund Ruin Probability

The aim is to study the fund ruin probability, i.e. the game reserves exhaustion probability. We consider $x$ and $k$ real numbers fulfilling $x > 0$ and $k > x$. First, the evaluation of $\rho_k(x)$, the probability that the fund reserves decrease from an initial value $x$ to a value in $(-\infty, 0]$ before reaching a value in $[k, +\infty)$, is considered. Then, after the calculus of the limit, as seen in the previous section, the evaluation of $\rho(x)$, the eventual fund ruin probability is considered, by admitting in such a situation that the random walk - that represents its reserves - evolves with no restrictions at the right of 0.

This method that was presented here is known in the literature of stochastic processes as Wald's Approximation. The explanations by Grimmett and Stirzaker [7] and Cox and Miller [16], relating this issue are carefully considered in our work. See also Ferreira, Andrade, Filipe, and Coelho [17]. Here, we consider the $S_n = \check{S}_n - x$ process, i.e., the random walk $S_0 = 0, S_n = S_{n-1} + X_n, n = 1, 2, ...$, instead of $\check{S}_n$ process.

Accordingly, when we evaluate $\rho_k(x)$, it is the probability the process $S_n$ is visiting the set $(-\infty, -x]$ before visiting the set $[k - x, +\infty)$ that, in reality, is considered. And when we evaluate $\rho(x)$ it is only the probability that the process $S_n$ goes down from the initial value 0 till a level lesser or equal than $-x$ that is considered.

We begin considering now the non-null value $\theta$ for which the $X_1$ moments generator function assumes the value 1. It is assumed that such a $\theta$ exists, that is, $\theta$ satisfies

$$E\left[e^{\theta X_1}\right] = 1, \theta \neq 0 \qquad (3.2.1).$$

Defining the process as $M_n = e^{\theta S_n}, n = 0, 1, 2, ...$, it is now evident that $E[|M_n|] < \infty$ and also, after (3.2.1), $E[M_{n+1}|X_1, X_2, ..., X_n] = E\left[e^{\theta(S_n + X_{n+1})}|X_1, X_2, ..., X_n\right] = e^{\theta S_n} E\left[e^{\theta X_{n+1}}|X_1, X_2, ..., X_n\right] = M_n$.

Consequently, the process $M_n$ is a Martingale in what concerns the sequence of random variables $X_1, X_2, ...$. Consider now $N$, the $S_n$ first passage time to outside the interval $(-x, k-x)$,

$$N = \min\{n \geq 0 : S_n \leq -x \text{ or } S_n \geq k - x\}.$$

The random variable $N$ is a stopping time – or a Markov time – for which the following conditions are assured: $\begin{cases} E[N] < \infty \\ E[|M_{n+1} - M_n||X_1, X_2, ..., X_n] \leq 2e^{|\theta|a} \end{cases}$, for $n < N$ and $a = -x$ or $a = k - x$. We recommend giving a look at the work by Grimmett and Stirzaker [7] and Figueira [18] on this issue. Under these conditions, we can resort to the Martingales Stopping Time Theorem and so:

$$E[M_n] = E[M_0] = 1 \qquad (3.2.2).$$

Besides,

$$E[M_n] = E\left[e^{\theta S_N}|S_N \leq -x\right]P(S_N \leq -x) + E\left[e^{\theta S_N}|S_N \geq k - x\right]P(S_N \geq k - x) \quad (3.2.3).$$

Realizing the approximations $E\left[e^{\theta S_N}|S_N \leq -x\right] \cong e^{-\theta x}$ and $E\left[e^{\theta S_N}|S_N \geq k - x\right] \cong e^{\theta(k-x)}$, and

considering $P(S_N \leq -x) = \rho_k(x) = 1 - P(S_N \geq k - x)$, after (3.2.2) and (3.2.3), we obtain

$$\rho_k(x) \cong \frac{1 - e^{\theta(k-x)}}{e^{-\theta x} - e^{\theta(k-x)}}, \text{when } E[X_1] \neq 0 \quad (3.2.4).$$

This result is the Classic Approximation for the Ruin Probability in the conditions stated in (3.2.1). To admit a non-null solution $\theta$ for the equation $E[e^{\theta X_1}] = 1$ implies indeed to assume that $E[X_1] \neq 0$.

Going farer, we may consider a particular case, beyond the studied above, looking at the situation for which the equation $E[e^{\theta X_1}] = 1$ only solution is precisely $\theta = 0$; it means, the situation at which $E[X_1] = 0$. This case may be considered through the following passage to the limit:

$$p_k(x) \cong \lim_{\theta \to 0} \frac{1 - e^{\theta(k-x)}}{e^{-\theta x} - e^{\theta(k-x)}} = \frac{k-x}{k}, \quad \text{when } E[X_1] = 0 \quad (3.2.5).$$

As for $p(x)$, the probability that the process $S_n$ decreases eventually from the initial value 0 to a level lesser or equal than $-x$, is also got from (3.2.4), now for a different passage to the limit:

$$p(x) \cong \lim_{k \to \infty} \frac{1 - e^{\theta(k-x)}}{e^{-\theta x} - e^{\theta(k-x)}} = e^{\theta x}, \quad \text{if } \theta < 0 \Leftrightarrow E[X_1] > 0 \quad (3.2.6).$$

Considering the previous section results on the simple random walk, it is effective to accept $p(x) = 1$ when $\theta \geq 0 \Leftrightarrow E[X_1] \leq 0$.

*3.3. A Ruin Probability's Particularization*

Consider $X_1, X_2, \ldots$ is a sequence of independent random variables with normal distribution with mean $\mu$ and standard deviation $\sigma$. So we can suppose $X_n$, is the fund reserves variation at the $n^{\text{th}}$ time unit, normally distributed with those parameters. Now, the moments' generator function is:

$$E[e^{\theta X_1}] = \frac{1}{\sqrt{2\pi}\sigma} \int_{-\infty}^{+\infty} e^{\theta x - \frac{(x-\mu)^2}{2\sigma^2}} dx = e^{\theta \mu + \frac{\theta^2 \sigma^2}{2}}. \text{ And equation (3.2.1) solution is } \theta = \frac{-2\mu}{\sigma^2}, \mu \neq 0.$$

So $p_k(x)$, the ruin probability, is acquired when substituting this result in (3.2.4), as follows:

$$\rho_k(x) \cong \frac{1 - e^{-\frac{2\mu(k-x)}{\sigma^2}}}{e^{\frac{2\mu x}{\sigma^2}} - e^{-\frac{2\mu(k-x)}{\sigma^2}}}, \text{ when } \mu \neq 0 \quad (3.3.1).$$

It is obvious that this particularization does not influence the approximation to $\rho_k(x)$ when $\mu = 0$. As it was seen before, it is given by (3.2.5) and after (3.2.6), we have

$$\rho(x) \cong e^{-\frac{2\mu x}{\sigma^2}}, \quad \text{when } \mu > 0 \quad (3.3.2).$$

*3.4. Assets* and *Liability Management Politics*

It is evident that management politics for the pensions fund must be defined to guarantee its sustainability.

Consider, for instance, a pensions fund management scheme that raises the assets value by some positive constant $\theta_n$, when the assets value falls equal to the liabilities process by the $n^{th}$ time. This corresponds to consider the modification $\bar{A}(t)$ of the process A(t) that restarts at times $T_n$ when A(t) hits the barrier L(t) by the $n^{th}$ time at the level $L(T_n) + \theta_n$.

It is a convenient choice the management policy where

$\theta_n = L(T_n)(e^\theta - 1)$, for some $\theta > 0$ (3.4.1).

The corresponding modification $\tilde{Y}(t)$ of Y(t) will behave as a generalized Brownian motion process that restarts at the level $\ln \frac{L(T_n) + \theta_n}{L(T_n)} = \theta$ when it hits 0 (at times $T_n$).

It is considered the pensions fund perpetual maintenance cost present value, because of the proposed asset-liability management scheme, given by the random variable:

$V(r, a, \theta) = \sum_{n=1}^{\infty} \theta_n e^{-rT_n} = \sum_{n=1}^{\infty} b(e^\theta - 1) e^{-(r-\rho)T_n}, r > \rho$, where r represents the appropriate discount interest rate. To obtain the above expression it was only made use of the L(t) definition and (3.4.1). Note that it is possible to express the expected value of the above random variable as

$$v_r(a, \theta) = \frac{b(e^\theta - 1)}{\theta} v_{r-\rho}(a, \theta) = \frac{b(e^\theta - 1) e^{-K_{r-\rho} a}}{1 - e^{-K_{r-\rho} \theta}} \quad (3.4.2).$$

As $\theta \to 0$

$$\lim_{\theta \to 0} v_r(a, \theta) = \frac{b e^{-K_{r-\rho} a}}{K_{r-\rho}} \quad (3.4.3).$$

In a similar way, the maintenance cost up to time t in the above-mentioned management scheme, is the stochastic process $W(t; r, a, \theta) = \sum_{n=1}^{N(t)} b(e^\theta - 1)e^{-(r-\rho)T_n}$, $W(t; r, a, \theta) = 0$ if $N(t) = 0$, with expected value function

$$w(t; a, \theta) = \frac{b(e^\theta - 1)}{\theta} w_{r-\rho}(t; a, \theta) \qquad (3.4.4).$$

## 4. Discussion and conclusions

Pension funds are important grounds of peoples' income during retirement. Once pension systems are important parts in financial markets as institutional investors, being the financial markets' participants aware of the possibility of having losses and gains, and once pensions funds need to keep viable standards, it is necessary and imperative to assure that pensions funds find a sustainability support, either by their investments in the market, by external injection of capital or by the reinforcement of internal contributions. That is why our theoretical methodology represents an important contribution for pensions funds management. The knowledge of the circumstances on which funds schemes develop and the contexts that they must face, with constraints and obstacles, and the possibility to go to bankruptcy imply that managers be well prepared to face funds reality troubles and crises. We show that it is essential to keep the long-term equilibrium of contributions and expenditures of a pensions scheme to achieve a sustainable system's evolution, as also expressed in several studies, as in Holzmann et al. [19]; Grech [20]; or Tian and Zhao [21].

It is true that pensions funds vary greatly from country to country, but fundamentally the basis is the same varying from the defined benefit to the defined contribution. Anyway, the basic principle is the same: the fund has to be sustainable and has to guarantee a future income to contributors. That is why a rigorous management process has to be maintained along the pensions fund life. This involves a meticulous definition of the criteria of funds' investments. Our study allows to understand how the fund is structured and how it may develop according to the probabilities based on the assumptions of the model.

Our study allows to get important conclusions, once the sustainability of a fund is vital for guaranteeing that older people, who are intended to be protected, can have their lives secured by the fund's investments and an assured economic independence when ageing.

Our paper, by studying the gambler's ruin probability (i.e. the probability of the game reserves exhaustion), gets the solution of the differential equation underlying a standard gambler game with a random walk assumption. This solution explains the suitability of pension fund systems. Also, the Brownian motion process, allowing to study how the fund evolves and its needs to be managed, gives results on the way and in what extent an intervention is required, and computing its costs.

The results in the study were obtained considering the simple and general random walk, and geometric and general Brownian motion, that are classic and widely studied stochastic processes. Since its general ideas are easily grasped by everyone, which quickly connect them with real systems, the random walk is used to model situations considering more disparate realities, far beyond the reserve evolution models considered in this work. They are also, and very often, used to build other complex systems, sometimes much more complex, for other types of systems.

We highlight in our approach a set of different methodologies that were applied to the study of this type of processes such as the cases of Difference Equations and Martingales Theory.

In our approach, reserves systems are treated as physical systems. In our study, we recognize that this may be a limitation to be considered, since it is not obvious that the direct application of these principles to financial reserve funds can be legitimate, when their own dynamics of appreciation and devaluation over time are ignored.

The models themselves, and the consequent valuation of stability systems based on the assessment of the probability of depletion or ruin of reserves, seem to be valid only in constant prices contexts.

In studying the managing scheme to try to guarantee the pensions fund stability, the interest rate, a factor associated with the process of temporal depreciation of the value of money when considering the modeling of financial reserves was considered. Analogous approaches can be seen, for instance, in Ferreira [11, 12]. Also, the demographic context that greatly influences the sustainability of retirement pensions must be considered, see, for instance Figueira and Ferreira [22].

Considering the results achieved, it seems evident that official authorities must consider the factors that influence the possible fluctuations on the pensions fund evolution and, resulting from there, they must plan the correct management in each period, and the proper strategy for that moment, being it of contracting economic activities' time or of a period of economic expansion. Our study helps to calculate the long-term costs of the various policies available to ensure the long-term sustainability of pension funds. This cost is the value of the fund's reinforcements, either through external capital injections or through the reinforcement of taxpayers' contributions, which are also the beneficiaries of the pension fund system. Actuarial reporting official systems must be operational for responding management's needs and, when required, solving pensions crises to come from difficult crises times and mostly from upcoming more aged populations.

**Author Contributions:** Conceptualization, M.F. and J.F.; methodology, M.F. and J.F.; software, M.F. and J.F.; validation, M.F. and J.F.; formal analysis, M.F. and J.F.; investigation, M.F. and J.F.; resources, M.F. and J.F..; data curation, M.F. and J.F.; writing—original draft preparation, M.F. and J.F.; writing—review and editing, M.F. and J.F.; visualization, M.F. and J.F.; supervision, M.F. and J.F.; project administration, M.F. and J.F.; funding acquisition, M.F. and J.F. All authors have read and agreed to the published version of the manuscript.

**Acknowledgments:** This work is financially supported by the Information Sciences, Technologies and Architecture Research Center (ISTAR-IUL), to which the authors are grateful.

**Conflicts of Interest:** The authors declare no conflict of interest.